\documentclass[12pt]{amsart}

\usepackage[margin=1.15in]{geometry}

 \usepackage{amssymb}
 \usepackage{graphicx}

\usepackage{amsthm,amssymb,amsmath,amscd,breqn}
\usepackage{amsfonts,color}
\usepackage{latexsym}
\usepackage[all,cmtip]{xy}
\usepackage{ulem}
\usepackage[colorlinks=true, pdfstartview=FitV, linkcolor=blue, citecolor=blue, urlcolor=blue]{hyperref}
\usepackage{marginnote}

\usepackage{mathptmx}

\newtheorem{thm}{Theorem}[section]
\newtheorem{lem}[thm]{Lemma}
\newtheorem{cor}[thm]{Corollary}
\newtheorem{prop}[thm]{Proposition}
\theoremstyle{definition}
\newtheorem{df}[thm]{Definition}
\theoremstyle{remark}
\newtheorem{rem}[thm]{Remark}

\newtheorem{example}[thm]{Example}

\newcommand{\bC}{{\mathbb C}}

\newcommand{\bH}{{\mathbb H}}

\newcommand{\bQ}{{\mathbb Q}}
\newcommand{\bZ}{{\mathbb Z}}

\newcommand{\cD}{{\mathcal D}}

\newcommand{\cF}{{\mathcal F}}

\newcommand{\cH}{{\mathcal H}}

\newcommand{\cL}{{\mathcal L}}

\newcommand{\cZ}{{\mathcal Z}}
\newcommand{\cP}{{\mathcal P}}
\newcommand{\cR}{{\mathcal R}}

\newcommand{\cQ}{{\mathcal Q}}

\newcommand{\cV}{{\mathcal V}}

\newcommand{\rk}{\mathop{{\mathrm{rank}}}\nolimits}

\newcommand{\lra}{\longrightarrow}

\newcommand{\cok}{\mathop{{\mathrm{coker}}}\nolimits}

\def\be{\begin{equation}}
\def\ee{\end{equation}}
\def\bt{\begin{thm}}
\def\et{\end{thm}}
\def\bc{\begin{cor}}
\def\ec{\end{cor}}
\def\br{\begin{rem}}
\def\er{\end{rem}}
\def\bp{\begin{prop}}
\def\ep{\end{prop}}
\def\bl{\begin{lem}}
\def\el{\end{lem}}
\def\bn{\begin{enumerate}}
\def\en{\end{enumerate}}
\def\bex{\begin{example}}
\def\eex{\end{example}}
\def\bd{\begin{df}}
\def\ed{\end{df}}

\begin{document}                        


\title[Vanishing cohomology]{Vanishing cohomology and Betti bounds for complex projective hypersurfaces}


\author[L. Maxim ]{Lauren\c{t}iu Maxim}
\address{L. Maxim : Department of Mathematics, University of Wisconsin-Madison, 480 Lincoln Drive, Madison WI 53706-1388, USA}
\email {maxim@math.wisc.edu}
\thanks{L. Maxim is partially supported by the Simons Foundation Collaboration (Grant \#567077) and by 
the Romanian Ministry of National Education (CNCS-UEFISCDI grant PN-III-P4-ID-PCE-2020-0029).}

\author[L. P\u{a}unescu ]{Lauren\c{t}iu P\u{a}unescu}
\address{L. P\u{a}unescu: Department of Mathematics, University of Sydney, Sydney, NSW, 2006, Australia}
\email {laurentiu.paunescu@sydney.edu.au}

\author[M. Tib\u{a}r]{Mihai Tib\u{a}r}
\address{M. Tib\u{a}r : Universit\' e de  Lille, CNRS, UMR 8524 -- Laboratoire Paul Painlev\'e, F-59000 Lille, France}  
\email {mihai-marius.tibar@univ-lille.fr}
\thanks{M. Tib\u{a}r acknowledges the support of the Labex CEMPI (ANR-11-LABX-0007). }

\keywords{singular projective hypersurface, vanishing cycles, vanishing cohomology, Betti numbers, Milnor fiber, Lefschetz hyperplane theorem}

\subjclass[2010]{32S30, 32S50, 55R55, 58K60}

\date{\today}

\begin{abstract} We employ the formalism of vanishing cycles and perverse sheaves to introduce and study the vanishing cohomology of complex projective hypersurfaces. As a consequence, we give upper bounds for the Betti numbers of projective hypersurfaces, generalizing those obtained by different methods by Dimca in the isolated singularities case, and by Siersma-Tib\u{a}r in the case of hypersurfaces with a $1$-dimensional singular locus. We also prove a supplement to the Lefschetz hyperplane theorem for hypersurfaces, which takes the dimension of the singular locus into account, and we use it to give a new proof of a result of Kato.
\end{abstract}


\maketitle   




\section{Introduction. Results}\label{intro}
Let $V=\{f=0\} \subset \bC P^{n+1}$ be a reduced complex projective hypersurface of degree $d$, with $n \geq 1$. By the classical Lefschetz Theorem, the inclusion map $j:V \hookrightarrow  \bC P^{n+1}$ induces cohomology isomorphisms 
\be\label{one}
j^k:H^k( \bC P^{n+1};\bZ) \overset{\cong}\lra H^k(V;\bZ) \ \  \text{for all} \ \  k<n,
\ee
 and a primitive monomorphism for $k=n$ (e.g., see \cite[Theorem 5.2.6]{Di}). Moreover, if $s=\dim V_{\rm sing}<n$ is the complex dimension of the singular locus of $V$ (with $\dim \emptyset=-1$), then Kato \cite{Ka} showed that (see also \cite[Theorem 5.2.11]{Di})
\be\label{two}
H^k(V;\bZ) \cong H^k( \bC P^{n+1};\bZ) \ \  \text{for all} \ \  n+s+2\leq k\leq 2n,
\ee and the homomorphism $j^k$ induced by inclusion is given in this range (and for $k$ even) by multiplication by $d=\deg(V)$. 
It therefore remains to study the cohomology groups $H^k(V;\bZ)$ for $n \leq k \leq n+s+1$.

In the case when $V\subset \bC P^{n+1}$ is a {\it smooth} degree $d$ hypersurface, the above discussion yields that $H^k(V;\bZ) \cong H^k( \bC P^{n};\bZ)$ for all $k \neq n$. This is in fact the only information we take as an input in this note (it also suffices to work with \eqref{one}, its homology counterpart, and Poincar\'e duality). 
The Universal Coefficient Theorem also yields in this case that  $H^n(V;\bZ)$ is free abelian, and its rank $b_n(V)$ can be easily deduced from the formula for the Euler characteristic of $V$ (e.g., see \cite[Proposition 10.4.1]{M}): 
\be\label{chi}
\chi(V)=(n+2)-\frac{1}{d} \big[1+(-1)^{n+1}(d-1)^{n+2}\big].
\ee
Specifically, if $V\subset \bC P^{n+1}$ is a smooth degree $d$ projective hypersurface, one has:
\be\label{bsm}
b_n(V)=\frac{(d-1)^{n+2}+(-1)^{n+1}}{d}+\frac{3(-1)^n+1}{2}.
\ee 

The case when $V$ has only isolated singularities was studied by Dimca \cite{Di0,Di}, (see also \cite{Mi} and \cite{ST}) while projective hypersurfaces with a one-dimensional singular locus have been more recently considered by Siersma-Tib\u{a}r \cite{ST}. 

In the singular case, let us fix a Whitney stratification $\cV$ of $V$ and consider a one-parameter smoothing of degree $d$, namely $$V_t:=\{f_t=f-tg=0\}\subset \bC P^{n+1}  \ \ (t \in \bC),$$ for $g$ a general polynomial of degree $d$. Here,  the meaning of ``general'' is that the hypersurface $W:=\{g=0\}$ is smooth and transverse to all strata in the stratification $\cV$ of $V$.
Then, for $t \neq 0$ small enough, all the $V_t$ are smooth and transverse to the stratification $\cV$. Let $$B=\{f=g=0\}$$ be the base locus (axis) of the pencil. Consider the incidence variety
$$V_D:=\{(x,t)\in \bC P^{n+1} \times D \mid x \in V_t \}$$
with $D$ a small disc centered at $0 \in \bC$ so that $V_t$ is smooth for all $t \in D^*:=D\setminus \{0\}$. Denote by $\pi:V_D \to D$ the proper projection map, and note that $V=V_0=\pi^{-1}(0)$ and $V_t=\pi^{-1}(t)$ for all $t \in D^*$. 
In what follows we write $V$ for $V_0$ and use $V_t$ for a smoothing of $V$ (i.e., with $t \in D^*$). 
In this setup, one can define the Deligne vanishing cycle complex of the family $\pi$, see \cite[Section 10.3]{M} for a quick introduction. More precisely, one has a bounded constructible complex $$\varphi_\pi \underline{\bZ}_{V_D} \in D^b_c(V)$$
on the hypersurface $V$, whose hypercohomology groups fit into a long exact sequence (called the {\it specialization sequence}):
\be\label{spec}
\cdots \lra H^k(V;\bZ) \overset{sp^k}{\lra} H^k(V_t;\bZ) \overset{\alpha^k}{\lra} \bH^k(V; \varphi_\pi \underline{\bZ}_{V_D}) \lra H^{k+1}(V;\bZ) \overset{sp^{k+1}}{\lra} \cdots
\ee
The maps $sp^k$ are called the {\it specialization} morphisms, while the $\alpha^k$'s are usually referred to as the {\it canonical} maps. 
For any integer $k$, we define $$H^k_\varphi(V):=\bH^k(V; \varphi_\pi \underline{\bZ}_{V_D})$$
and call it the {\it $k$-th vanishing cohomology group of $V$}. This is an invariant of $V$, i.e., it does not depend on the choice of a particular smoothing of degree $d$ (since all smooth hypersurfaces of a fixed degree are diffeomorphic). By its very definition, the vanishing cohomology measures the difference between the topology of a given projective hypersurface $V$ and that of a smooth hypersurface of the same degree.

\br\label{r1}
Since the incidence variety $V_D=\pi^{-1}(D)$ deformation retracts to $V=\pi^{-1}(0)$, and the {\it specialization map}  $sp^k:H^k(V;\bZ) \to H^k(V_t;\bZ)$ of \eqref{spec} factorizes as $$H^k(V;\bZ) \overset{\cong}{\lra} H^k(V_D;\bZ) \lra H^k(V_t;\bZ)$$ with $H^k(V_D;\bZ) \to H^k(V_t;\bZ)$ induced by the inclusion map, it follows readily that the vanishing cohomology of $V$ can be identified with the relative cohomology of the pair $(V_D,V)$, i.e., \be H^k_\varphi(V) \cong H^{k+1}(V_D,V_t;\bZ).\ee
In particular, the groups $H^k_\varphi(V)$ are the cohomological version of the {\it vanishing homology groups} $$H_k^{\curlyvee}(V):=H_{k}(V_D,V_t;\bZ)$$ introduced and studied in \cite{ST} in special situations. For the purpose of computing Betti numbers of projective hypersurfaces, the two ``vanishing'' theories yield the same answer, but additional care is needed to handle torsion when computing the actual integral cohomology groups.
\er

Our first result gives the concentration degrees of the vanishing cohomology of a projective hypersurface in terms of the dimension of the singular locus. 
\bt\label{th1}
Let $V \subset \bC P^{n+1}$ be a reduced complex projective hypersurface with $s=\dim V_{\rm sing}$ the complex dimension of its singular locus. Then
\be
H^k_\varphi(V) \cong 0 \ \  \text{ for all integers} \ \ k \notin [n, n+s].
\ee
Moreover, $H^n_\varphi(V)$ is a free abelian group.
\et

In view of Remark \ref{r1}, one gets by Theorem \ref{th1} and the Universal Coefficient Theorem the concentration degrees of the vanishing homology groups $H_k^{\curlyvee}(V)$ of a projective hypersurface in terms of the dimension of its singular locus:
\bc\label{corh}
With the above notations and assumptions, we have that
\be
H_k^{\curlyvee}(V) \cong 0 \ \  \text{ for all integers} \ \ k \notin [n+1, n+s+1].
\ee
Moreover, $H_{n+s+1}^{\curlyvee}(V)$ is free.
\ec

\br
In the case when the projective hypersurface $V \subset \bC P^{n+1}$ has a $1$-dimensional singular locus, it was shown in \cite[Theorem 4.1]{ST} that $H_k^{\curlyvee}(V) \cong 0$ for all $k \neq n+1, n+2$. Moreover, Theorem 6.1 of \cite{ST} shows that in this case one also has that $H_{n+2}^{\curlyvee}(V)$ is free. So, Corollary \ref{corh} provides a generalization of the results of \cite{ST} to projective hypersurfaces with arbitrary singularities. Nevertheless, the methods used in its proof  are fundamentally different from those in \cite{ST}. 
\er

As a consequence of Theorem \ref{th1}, the specialization sequence \eqref{spec} together with the fact that the integral cohomology of a smooth projective hypersurface is free, yield the following result on the integral cohomology of a complex projective hypersurface (where the estimate on the $n$-th Betti number uses formula \eqref{bsm}):
\bc\label{corgen}
Let $V \subset \bC P^{n+1}$ be a degree $d$ reduced projective hypersurface with a singular locus $V_{\rm sing}$ of complex dimension $s$. Then:
\begin{itemize}
\item[(i)] $H^k(V;\bZ) \cong H^k(V_t;\bZ) \cong H^k(\bC P^n;\bZ)$ \ for all integers $k \notin [n, n+s+1]$.
\item[(ii)] $H^n (V;\bZ) \cong \ker(\alpha^n)$ is free.
\item[(iii)] $H^{n+s+1} (V;\bZ)  \cong H^{n+s+1} (\bC P^n;\bZ) \oplus \cok(\alpha^{n+s})$.
\item[(iv)] $H^k(V;\bZ) \cong \ker(\alpha^k) \oplus \cok(\alpha^{k-1})$ for all integers $k \in [n+1,n+s]$, $s\ge 1$.

\end{itemize}
In particular, $$b_n(V) \leq b_n(V_t)=\frac{(d-1)^{n+2}+(-1)^{n+1}}{d}+\frac{3(-1)^n+1}{2},$$
and 
$$b_{k}(V) \leq \rk \ H^{k-1}_{\varphi}(V) +  b_{k}(\bC P^{n})  \ \ \text{ for all integers} \ k \in [n+1,n+s+1],  \ s\ge 0.$$
\ec

The homological version of the  specialisation sequence \eqref{spec} identifies to the long exact sequence of the pair $(V_{D}, V_{t})$, namely: 
\[     \cdots \to H_{k+1}(V_{t};\bZ) \to H_{k+1}(V_{D};\bZ) \to H_{k+1}^{\curlyvee}(V;\bZ)  \ \overset{\alpha_{k}}{\longrightarrow} \  H_{k}(V_{t};\bZ) \to \cdots
\]
The  inclusions $V_{t}\hookrightarrow V_{D}\hookrightarrow \bC P^{n+1}\times D$ induce in homology a commutative triangle, where $H_{k}(V_{t};\bZ) \to H_{k}(\bC P^{n+1}\times D;\bZ)$ is injective for  $k\not= n$ (by the Lefschetz Theorem for $k<n$, and
it is multiplication by $d$ for $k>n$, see e.g. Remark \ref{Katoh}  for the homological version of the proof of Theorem \ref{Kato}). This
shows that the morphisms $H_{k}(V_{t};\bZ) \to H_{k}(V_{D};\bZ)$ is also injective for all $k\not= n$, and therefore $\alpha_{k} =0$ for $k\not= n$. Consequently,  the above long exact sequence splits into a $5$-term exact sequence,  and short exact sequences:
\begin{equation}\label{sp1}
\begin{split}
 0 \to  H_{n+1}(V_{t};\bZ)  \to  H_{n+1}(V;\bZ)    \to  H_{n+1}^{\curlyvee}(V;\bZ)    \stackrel{\alpha_{n}}{\to} 
 H_{n}(V_{t};\bZ)  \rightarrow H_{n}(V;\bZ)  \to 0. \\
 0 \to  H_{k}(V_{t};\bZ) \to H_{k}(V_{D};\bZ) \to H_{k}^{\curlyvee}(V;\bZ)  \to 0
\ \ \ \ \ \ \ \ \mbox{ for } \ k\ge n+1.
\end{split}
\end{equation}
We then get the following homological version of Corollary \ref{corgen}(i-iv), with the same upper bounds for Betti numbers, but with an interesting improvement for (iii) and (iv) showing more explicitly the dependence of the homology of $V$ on the vanishing homology groups:

\bc\label{corgenhom} 
Let $V \subset \bC P^{n+1}$ be a degree $d$ reduced projective hypersurface with a singular locus $V_{\rm sing}$ of complex dimension $s$. Then:
 \begin{itemize}
\item[(i')] $H_{k}(V;\bZ) \cong H_{k}(V_t;\bZ) \cong H_{k}(\bC P^n;\bZ)$ \ for all  $k\le n-1$ and  all $k\ge n+s+2$.
\item[(ii')] $H_{n}(V;\bZ) \cong \cok(\alpha_{n})$.
\item[(iii')] $H_{n+1}(V;\bZ) \cong \ker(\alpha_{n})  \oplus H_{n+1}(\bC P^n;\bZ)$.
\item[(iv')] $H_{k}(V;\bZ) \cong H_k^{\curlyvee}(V;\bZ) \oplus H_{k}(\bC P^n;\bZ)$, for all $ n+2 \le k \le n+s+1$, whenever $s\ge 1$, \\ and  $H_{n+s+1}(V;\bZ)$ is free.
\end{itemize} 
\ec

The ranks of the (possibly non-trivial) vanishing (co)homology groups can be estimated in terms of the local topology of singular strata and of their generic transversal types by making use of the hypercohomology spectral sequence. Such estimates can be made precise for hypersurfaces with low-dimensional singular loci. 
 Concretely, as special cases of Corollaries \ref{corgen} and \ref{corgenhom}, in Section \ref{bounds}  we recast Siersma-Tib\u{a}r's  \cite{ST} result for $s\le 1$, and in particular Dimca's  \cite{Di0,Di} computation  for $s=0$.  Concerning the estimation of the rank of the highest interesting (co)homology group, we prove the following general result:

\bt\label{th2}
Let $V \subset \bC P^{n+1}$ be a degree $d$ reduced projective hypersurface with a singular locus $V_{\rm sing}$ of complex dimension $s$. For each connected stratum $S_i \subseteq  V_{\rm sing}$ of top dimension $s$ in a Whitney stratification of $V$, let $F_i^\pitchfork$ denote its transversal Milnor fiber with corresponding Milnor number $\mu_i^\pitchfork$.
Then:
\be\label{bt}
b_{n+s+1}(V) \leq 1+ \sum_i \mu_i^\pitchfork,
\ee
 and the  inequality is strict for $n+s$ even.
\et

 In fact, the inequality in \eqref{bt} is deduced from  
 \be\label{btb} b_{n+s+1}(V) \leq 1+\rk \ H^{n+s}_{\varphi}(V),\ee
together with
\be
\rk \ H^{n+s}_{\varphi}(V) \leq \sum_i \mu_i^\pitchfork,
\ee
and the inequality \eqref{btb} is strict for $n+s$ even. For further refinements of Theorem \ref{th2}, see Remark \ref{rem31}. Note also that if $s=0$, i.e., $V$ has only isolated singularities, then $\mu_i^\pitchfork$ is just the usual Milnor number of such a singularity of $V$.

Let us remark that if the projective hypersurface $V \subset \bC P^{n+1}$ has singularities in codimension $1$, i.e., $s=n-1$, then $b_{n+s+1}(V)=b_{2n}(V)=r$, where $r$ denotes the number of irreducible components of $V$. Indeed, in this case, one has (e.g., see \cite[(5.2.9)]{Di}):
\be\label{top} H^{2n}(V;\bZ) \cong \bZ^{r}.\ee
In particular, Theorem \ref{th2} yields the following generalization of \cite[Corollary 7.6]{ST}:
\bc If the reduced projective hypersurface $V \subset \bC P^{n+1}$ has singularities in codimension $1$, then the number $r$ of irreducible components of $V$ satisfies the inequality:
\be
r \leq 1+\sum_i \mu_i^\pitchfork.
\ee
\ec

\br\label{fr} Note that if the projective hypersurface $V \subset \bC P^{n+1}$ is a rational homology manifold, then the Lefschetz isomorphism \eqref{one} and Poincar\'e duality over the rationals yield that $b_i(V)=b_i(\bC P^n)$ for all $i \neq n$. Moreover, $b_n(V)$ can be deduced by computing the Euler characteristic of $V$, e.g., as in \cite[Section 10.4]{M}.\er

The computation of Betti numbers  of a projective hypersurface which is a rational homology manifold can be deduced without appealing to Poincar\'e duality by using the vanishing cohomology instead, as the next result shows:
\bp\label{p1}
If the projective hypersurface $V \subset \bC P^{n+1}$ is a $\bQ$-homology manifold, then $H^k_{\varphi}(V) \otimes \bQ \cong 0$ for all $k \neq n$. In particular, in this case one gets: $b_i(V)=b_i(V_t)=b_i(\bC P^n)$ for all $i \neq n$, and $b_n(V)=b_n(V_t)+\rk H^n_{\varphi}(V)$.
\ep

\medskip

At this end, we note that Corollary \ref{corgen}(i) reproves Kato's isomorphism \eqref{two} about the integral cohomology of $V$, by using only the integral cohomology of a smooth hypersurface (for this it suffices to rely only on the Lefschetz isomorphism \eqref{one}, its homological version, and Poincar\'e duality). In Section \ref{supLHT}, we give a new proof of Kato's result (see Theorem \ref{Kato}), 
which relies on the following supplement to the Lefschetz hyperplane section theorem for hypersurfaces, which may be of independent interest: 
\bt\label{thapi}
Let $V \subset \bC P^{n+1}$ be a reduced complex projective hypersurface with $s=\dim V_{\rm sing}$ the complex dimension of its singular locus. (By convention, we set $s=-1$ if $V$ is nonsingular.) Let $H \subset \bC P^{n+1}$ be a generic hyperplane. 
Then
\be\label{34api}
H^k(V,V\cap H; \bZ)=0 \ \ \text{for} \ \ k < n \ \ \text{and} \ \ n+s+1 < k < 2n.
\ee
Moreover, $H^{2n}(V,V \cap H; \bZ)\cong\bZ^r$, where $r$ is the number of irreducible components of $V$, and $H^{n}(V,V \cap H; \bZ)$ is (torsion-)free.
\et

Note that the vanishing \eqref{34api} for $k<n$ is equivalent to the classical Lefschetz hyperplane section theorem. The proof of \eqref{34api} for $n+s+1 < k < 2n$ reduces to understanding the homotopy type of the complement of a smooth affine hypersurface transversal to the hyperplane at infinity;
see \cite[Corollary 1.2]{Lib}  for such a description. Homological counterparts of Theorem \ref{thapi} and of Kato's result are also explained in Section \ref{supLHT}, see Corollary \ref{corap} and Remark \ref{Katoh}.

\medskip

Finally, let us note that similar techniques apply to the study of Milnor fiber cohomology of complex hypersurface singularity germs.  This is addressed by the authors in the follow-up paper \cite{MPT} (see also \cite{ST0} for the case of $1$-dimensional singularities).

\medskip

\noindent{\bf Acknowledgements.} L. Maxim thanks the Sydney Mathematical Research Institute (SMRI) for support and hospitality, and J\"org Sch\"urmann for useful discussions. 

\section{Concentration degrees of vanishing cohomology}

The proof of Theorem \ref{th1} makes use of the formalism of perverse sheaves and their relation to vanishing cycles, see \cite{Di1,M} for a brief introduction.

\subsection{Proof of Theorem \ref{th1}}
By definition, the incidence variety $V_D$ is a complete intersection of pure complex dimension $n+1$. It is non-singular if $V=V_0$ has only isolated singularities, but otherwise it has singularities where the base locus $B=V\cap W$ of the pencil  $\{f_t\}_{t\in D}$ intersects the singular locus $\Sigma:=V_{\rm sing}$ of $V$.

If $\underline{\bZ}_{V_D}$ denotes the constant sheaf with stalk $\bZ$ on the complete intersection $V_D$, a result of L\^e \cite{Le} implies that the complex $\underline{\bZ}_{V_D}[n+1]$ is a perverse sheaf on $V_D$. It then follows that 
$\varphi_\pi \underline{\bZ}_{V_D}[n]$ is a $\bZ$-perverse sheaf on $\pi^{-1}(0)=V$ (see, e.g., \cite[Theorem 10.3.13]{M} and the references therein).  

Recall that the stalk of the cohomology sheaves of $\varphi_\pi \underline{\bZ}_{V_D}$ at a point $x \in V$ are computed by (e.g., see \cite[(10.20)]{M}):
\be
\cH^j(\varphi_\pi \underline{\bZ}_{V_D})_x \cong H^{j+1}(B_{x}, B_{x}\cap V_t;\bZ),
\ee
where $B_{x}$ denotes the intersection of $V_D$ with a sufficiently small ball in some chosen affine chart $\bC^{n+1} \times D$ of the ambient space $\bC P^{n+1} \times D$ (hence $B_x$ is contractible). Here $B_{x}\cap V_t=F_{\pi,x}$ is the Milnor fiber of $\pi$ at $x$.
Let us now consider the function $$h=f/g:\bC P^{n+1} \setminus W \to \bC$$ where $W:=\{g=0\}$, and note that $h^{-1}(0)=V\setminus B$ with $B=V\cap W$  the base locus of the pencil. If $x \in V \setminus B$, then in a neighborhood of $x$ one can describe $V_t$  ($t \in D^*$)  as $$\{x \mid f_t(x)=0\}=\{x \mid h(x)=t\},$$ 
i.e., as the Milnor fiber of $h$ at $x$. Note also that $h$ defines $V$ in a neighborhood of $x \notin B$. Since the Milnor fiber of a complex hypersurface singularity germ does not depend on the choice of a local equation (e.g., see \cite[Remark 3.1.8]{Di}), we can therefore use $h$ or a local representative of $f$ when considering Milnor fibers (of $\pi$) at points in $V \setminus B$. From here on we will use the notation $F_x$ for the Milnor fiber of the hypersurface singularity germ $(V,x)$, and we note for future reference that the above discussion also yields that $F_x$ is a manifold, which moreover is contractible if $x \in V \setminus B$ is a smooth point.

It was shown in \cite[Proposition 5.1]{PP} (see also  \cite[Proposition 4.1]{MSS} or \cite[Lemma 4.2]{ST}) that
there are no vanishing cycles along the base locus $B$, i.e., 
\be \varphi_\pi \underline{\bZ}_{V_D} \vert_B \simeq 0.\ee
Therefore, if $u:V\setminus B \hookrightarrow V$ is the open inclusion, we get that 
\be\label{s6}
\varphi_\pi \underline{\bZ}_{V_D} \simeq u_! u^*  \varphi_\pi \underline{\bZ}_{V_D}.
\ee 
Since pullback to open subvarieties preserves perverse sheaves, we note that $u^*  \varphi_\pi \underline{\bZ}_{V_D}[n]$ is a perverse sheaf on the {\it affine} variety $V \setminus B$. 
Artin's vanishing theorem for perverse sheaves (e.g., \cite[Corollary 6.0.4]{Sc}) then implies that:

\begingroup
\allowdisplaybreaks
\begin{equation}\label{s1}
\begin{split}
 H^k_\varphi(V) 
 & := \bH^{k}(V; \varphi_\pi \underline{\bZ}_{V_D}) \\
  & \cong \bH^{k-n}(V; \varphi_\pi \underline{\bZ}_{V_D}[n]) \\
 & \cong \bH^{k-n}(V; u_! u^* \varphi_\pi \underline{\bZ}_{V_D}[n]) \\
 & \cong \bH_c^{k-n}(V \setminus B; u^* \varphi_\pi \underline{\bZ}_{V_D}[n]) \\
& \cong 0
\end{split}
\end{equation}
\endgroup
 for all $k-n<0$, or equivalently, for all $k<n$. 

Contractibility of Milnor fibers at smooth points of $V \setminus B$ implies that the support of $\varphi_\pi \underline{\bZ}_{V_D}$ is in fact contained in $\Sigma \setminus B$, with $\Sigma$ denoting as before the singular locus of $V$. In particular, if $v:\Sigma \setminus B \hookrightarrow V\setminus B$ is the closed inclusion, then 
\be\label{s5} u^* \varphi_\pi \underline{\bZ}_{V_D} \simeq v_!v^*u^* \varphi_\pi \underline{\bZ}_{V_D}.\ee
Next, consider the composition of inclusion maps $$\Sigma \setminus B \overset{q}{\hookrightarrow}\Sigma \overset{p}{\hookrightarrow} V$$ with $p\circ q=u \circ v$. By using \eqref{s6} and \eqref{s5}, we get:
\begingroup
\allowdisplaybreaks
\begin{equation}\label{s3}
\begin{split} \varphi_\pi \underline{\bZ}_{V_D} 
& \simeq u_!v_!v^* u^*\varphi_\pi \underline{\bZ}_{V_D} \\
& \simeq (u\circ v)_! (u\circ v)^* \varphi_\pi \underline{\bZ}_{V_D} \\
& \simeq (p\circ q)_! (p\circ q)^* \varphi_\pi \underline{\bZ}_{V_D} \\
& \simeq p_!q_!q^*p^* \varphi_\pi \underline{\bZ}_{V_D} \\
&\simeq p_*p^* \varphi_\pi \underline{\bZ}_{V_D},
\end{split}
\end{equation}
\endgroup
where the last isomorphism uses the fact that $p^* \varphi_\pi \underline{\bZ}_{V_D}$ is supported on $\Sigma \setminus B$, hence $p^* \varphi_\pi \underline{\bZ}_{V_D}\simeq q_!q^*p^* \varphi_\pi \underline{\bZ}_{V_D}$.
Since the support of the perverse sheaf $\varphi_\pi \underline{\bZ}_{V_D}[n]$ on $V$ is contained in the closed subset $\Sigma$, we get that $p^*\varphi_\pi \underline{\bZ}_{V_D}[n]$ is a perverse sheaf on $\Sigma$ (e.g., see \cite[Corollary 8.2.10]{M}). Since the complex dimension of $\Sigma$ is $s$, the support condition for perverse sheaves together with the hypercohomology spectral sequence yield that 
$$\bH^{\ell}(\Sigma; p^*\varphi_\pi \underline{\bZ}_{V_D}[n]) \cong 0$$ for all $\ell \notin [-s,s]$. 
This implies by \eqref{s3} that
\be\label{s2} H^k_\varphi(V)=\bH^{k-n}(V; \varphi_\pi \underline{\bZ}_{V_D}[n]) \cong 
\bH^{k-n}(\Sigma;p^*\varphi_\pi \underline{\bZ}_{V_D}[n]) \cong 0\ee for all $k \notin [n-s,n+s]$.

The desired concentration degrees for the vanishing cohomology is now obtained by combining \eqref{s1} and \eqref{s2}. 

\medskip

Let us finally show that $H^n_\varphi(V)$ is free. Fix a Whitney stratification $\cV$ of $V$, so that $V \setminus \Sigma$ is the top stratum. (Note that together with $\pi^{-1}(D^*)$, this also yields a Whitney stratification of $V_D$.) Since $W$ intersects $V$ transversally (i.e., $W$ intersects each stratum $S$ in $\cV$ transversally in $\bC P^{n+1}$), we can assume without any loss of generality that the base locus $B=V \cap W$ is a closed union of strata of $\cV$.
Next, we have by \eqref{s1} 
that
$$ H^n_\varphi(V) \cong \bH_c^{0}(V \setminus B; u^* \varphi_\pi \underline{\bZ}_{V_D}[n]),$$
with $$\cP:=u^* \varphi_\pi \underline{\bZ}_{V_D}[n]$$ 
a $\bZ$-perverse sheaf on the affine variety $V \setminus B$ and $u:V \setminus B \hookrightarrow V$ the open inclusion. In particular, this implies that if $S \in \cV$ is any stratum in $V \setminus B$ with inclusion $i_S:S \hookrightarrow V \setminus B$ then $\cH^k(i_S^!\cP) \simeq 0$ for all integers $k<-\dim_{\bC} S$.
By the Artin-Grothendieck type result of \cite[Corollary 6.0.4]{Sc}, in order to show that $\bH^0_c(V\setminus B;\cP)$ is free it suffices to check that the perverse sheaf $\cP$ satisfies the following costalk condition (see \cite[Example 6.0.2(3)]{Sc}):\footnote{We thank J\"org Sch\"urmann for indicating the relevant references to us.}
\be\label{co1}
\cH^{-\dim_{\bC} S}(i_S^!\cP)_x \ \text{ is free }
\ee
for any point $x$ in any stratum $S$ in $V \setminus B$ with inclusion $i_S:S \hookrightarrow V \setminus B$. 
Let us now fix a stratum $S \in \cV$ contained in $V \setminus B$ and let $x \in S$ be a point with inclusion map $k_x:\{x\} \hookrightarrow S$. Consider the composition $i_x:=i_S \circ k_x: \{x\} \hookrightarrow V \setminus B$. Using the fact that $$k_x^*i_S^! \simeq k_x^!i_S^! [2 \dim_{\bC} S] \simeq i_x^! [2 \dim_{\bC} S]$$
(e.g., see \cite[Remark 6.0.2(1)]{Sc}), the condition \eqref{co1} for $x \in S$ is equivalent to the following:
\be\label{co2}
\cH^{\dim_{\bC} S}(i_x^!\cP) \ \text{ is free}.
\ee
In fact, the above discussion applies to any algebraically constructible complex $\cF^\centerdot \in {^pD}^{\geq 0}$, with $({^pD}^{\leq 0}, {^pD}^{\geq 0})$ denoting the perverse t-structure on $D^b_c(V \setminus B)$.
Furthermore, in our setup (i.e., working with PID coefficients and having finitely generated stalk cohomology) $\cF^\centerdot \in {^pD}^{\geq 0}$ satisfies the additional costalk condition \eqref{co1} (or, equivalently, \eqref{co2}) 
if and only if the Verdier dual $\cD\cF^\centerdot $ satisfies $\cD\cF^\centerdot \in {^pD}^{\leq 0}$.

Let $i:V=V_0 \hookrightarrow V_D$ denote the closed inclusion, and consider the following {\it variation triangle} for the projection map $\pi:V_D \to D$:
\be\label{var}
i^![1] \lra \varphi_\pi \overset{var}{\lra} \psi_\pi \overset{[1]}{\lra}
\ee
with $\psi_\pi $ denoting the corresponding nearby cycle functor for $\pi$ (e.g., see \cite[(5.90)]{Sc}). Apply the functor $u^!=u^*$ to the triangle \eqref{var}, and the apply the resulting triangle of functors to the complex $\underline{\bZ}_{V_D}[n]$ to get the following triangle of constructible complexes on $V\setminus B$:
\be\label{var2}
\cZ:=u^!i^!\underline{\bZ}_{V_D}[n+1] \lra \cP:=u^* \varphi_\pi \underline{\bZ}_{V_D}[n] \lra 
\cR:=u^* \psi_\pi \underline{\bZ}_{V_D}[n] \overset{[1]}{\lra}
\ee
Let $x \in S$ be a point in a stratum of $V \setminus B$ with inclusion map $i_x:\{x\} \hookrightarrow V \setminus B$ as before, and apply the functor $i_x^!$ to the triangle \eqref{var2} to get the triangle:
\be\label{var3}
i_x^!\cZ \lra i_x^!\cP \lra i_x^!\cR \overset{[1]}{\lra}
\ee
The cohomology long exact sequence associated to \eqref{var3} contains the terms
$$\cdots \lra \cH^{\dim_{\bC} S}(i_x^!\cZ) \lra \cH^{\dim_{\bC} S}(i_x^!\cP) \lra \cH^{\dim_{\bC} S}(i_x^!\cR) \lra \cdots$$
Since the category of (torsion-)free abelian groups is closed under extensions, in order to prove \eqref{co2} it suffices to check that $\cH^{\dim_{\bC} S}(i_x^!\cZ)$ and $\cH^{\dim_{\bC} S}(i_x^!\cR)$ are (torsion-)free. (Note that, in fact, all costalks in question are finitely generated.)

Let us first show that $\cH^{\dim_{\bC} S}(i_x^!\cZ)$ is free. Regard the stratum $S$ containing $x$ as a stratum in $V_D$, and let $r_x:\{x\} \to V_D$ be the point inclusion, i.e., $r_x=i\circ u \circ i_x$. So $i_x^!\cZ=r_x^!\underline{\bZ}_{V_D}[n+1]$. 
Recall that $\underline{\bZ}_{V_D}[n+1]$ is a $\bZ$-perverse sheaf on $V_D$, i.e., $\underline{\bZ}_{V_D}[n+1] \in {^pD}^{\leq 0}(V_D)  \cap{^pD}^{\geq 0}(V_D)$.
As already indicated above, in order to show that $\cH^{\dim_{\bC} S}(r_x^!\underline{\bZ}_{V_D}[n+1])$ is free it suffices to verify that $\cD(\underline{\bZ}_{V_D}[n+1]) \in {^pD}^{\leq 0}(V_D)$, or equivalently, $\cD\underline{\bZ}_{V_D}  \in {^pD}^{\leq -n-1}(V_D)$. This fact is a consequence of \cite[Definition 6.0.4, Example 6.0.11]{Sc}, where it is shown that the complete interesection $V_D$ has a {\it rectified homological depth} equal to its complex dimension $n+1$.


Next note that, due to the local product structure, the Milnor fiber $F_x$ of the hypersurface singularity germ $(V,x)$ with $x \in S$ has the homotopy type of a finite CW complex of real dimension $n-\dim_{\bC} S$. In particular, $H_{n-\dim_{\bC}S}(F_x;\bZ)$ is free.
Since by the costalk calculation (cf. \cite[(5.92)]{Sc}) and Poincar\'e duality we have for $x\in S$ that
\be
\cH^{\dim_{\bC} S}(i_x^!\cR) \cong H_c^{n+\dim_{\bC}S}(F_x;\bZ) \cong H_{n-\dim_{\bC}S}(F_x;\bZ),
\ee
it follows that $\cH^{\dim_{\bC} S}(i_x^!\cR)$ is free. This completes the proof of Theorem \ref{th1}.

\subsection{Proof of Proposition \ref{p1}}
Since $V$ is a $\bQ$-homology manifold, it follows by standard arguments involving the Hamm fibration (e.g., see \cite[Theorem 3.2.12]{Di}) that $V_D$ is also a $\bQ$-homology manifold (with boundary). Thus $\underline{\bQ}_{V_D}[n+1]$ is a self-dual $\bQ$-perverse sheaf on $V_D$. Moreover, since $\varphi_\pi[-1]$ commutes with the Verdier dualizing functor (see \cite[Theorem 3.1]{Ma} and the references therein), we get that $\cQ:=\varphi_\pi\underline{\bQ}_{V_D}[n]$ is a Verdier self-dual perverse sheaf on $V$. Using the Universal Coefficients Theorem, we obtain:
$$H^k_{\varphi}(V) \otimes \bQ=\bH^{k-n}(V;\cQ)\cong \bH^{k-n}(V;\cD\cQ) \cong \bH^{n-k}(V;\cQ)^\vee = (H^{2n-k}_{\varphi}(V)\otimes \bQ)^\vee.$$
The desired vanishing follows now from Theorem \ref{th1}.


\section{Bounds on Betti numbers of projective hypersurfaces}\label{bounds}
In this section, we prove Theorem \ref{th2} and specialize it, along with Corollary \ref{corgen}, in the case when the complex dimension $s$ of the singular locus is $\leq 1$.

\subsection{Proof of Theorem \ref{th2}}
Let $\Sigma:=V_{\rm sing}$ be the singular locus of $V$, of complex dimension $s$, and fix a Whitney stratification $\cV$ of $V$ so that $V \setminus \Sigma$ is the top open stratum. 
We have by Corollary \ref{corgen} (or by the specialization sequence \eqref{spec}) that
$$b_{n+s+1}(V) \leq 1+\rk \ H^{n+s}_{\varphi}(V).$$
So it suffices to show that
\be\label{rvg}
\rk \ H^{n+s}_{\varphi}(V) \leq \sum_i \mu_i^\pitchfork,
\ee
where the summation on the right-hand side runs over the top $s$-dimensional connected strata $S_i$ of $\Sigma$, and $\mu_i^\pitchfork$ denotes the corresponding transversal Milnor number for such a stratum $S_i$.

If $s=0$, an easy computation shows that \eqref{rvg} is in fact an equality, see \eqref{nee} below. Let us next investigate the case when $s\geq 1$.

For any $\ell \leq s$, denote by $\Sigma_{\ell}$ the union of strata in $\Sigma$ of complex dimension $\leq \ell$. In particular, we can filter $\Sigma$ by closed (possibly empty) subsets
$$\Sigma=\Sigma_s \supset \Sigma_{s-1} \supset \cdots \supset \Sigma_0 \supset \Sigma_{-1}=\emptyset.$$
Let $$U_\ell:=\Sigma_{\ell} \setminus \Sigma_{\ell-1}$$ be the union of $\ell$-dimensional strata, so $\Sigma_\ell=\sqcup_{k\leq \ell} U_k$. (Here, $\sqcup$ denotes disjoint union.) Recall that the smooth hypersurface $W=\{g=0\}$ was chosen so that it intersects each stratum in $\Sigma$ transversally.

In the notations of the proof of Theorem \ref {th1} , it follows from equations \eqref{s1} and \eqref{s5} that:
$$ H^{n+s}_\varphi(V) \cong \bH_c^{n+s}(V \setminus B; u^* \varphi_\pi \underline{\bZ}_{V_D})
\cong \bH_c^{n+s}(\Sigma \setminus B; v^* u^*\varphi_\pi \underline{\bZ}_{V_D}),$$
with $B=V \cap W$ the axis of the pencil, and with $v:\Sigma\setminus B \hookrightarrow V \setminus B$ and $u:V \setminus B \hookrightarrow V$ the inclusion maps. 
We also noted that either $h$ or a local representative of $f$ can be used when considering Milnor fibers of $\pi$ at points in $V \setminus B$. For simplicity, let us use the notation
$$\cR := v^* u^*\varphi_\pi \underline{\bZ}_{V_D} \in D^b_c(\Sigma \setminus B),$$
and consider the part of the long exact sequence for the compactly supported hypercohomology of $\cR$ associated to the disjoint union 
$$\Sigma \setminus B= (U_s \setminus B) \sqcup (\Sigma_{s-1} \setminus B)$$
involving $H^{n+s}_{\varphi}(V)$, 
namely:
$$\cdots \to  \bH^{n+s}_c(U_s \setminus B; \cR) \to H^{n+s}_{\varphi}(V) \to \bH^{n+s}_c(\Sigma_{s-1} \setminus B; \cR) \to \cdots$$
  We claim that 
 \be\label{cl1}
 \bH^{n+s}_c(\Sigma_{s-1} \setminus B; \cR) \cong 0,
 \ee
 so, in particular, there is an epimorphism:
 \be\label{nee2}\bH^{n+s}_c(U_s \setminus B; \cR) \twoheadrightarrow H^{n+s}_{\varphi}(V).\ee

 In order to prove \eqref{cl1}, consider 
  the part of the long exact sequence for the compactly supported hypercohomology of $\cR$ associated to the disjoint union 
$$\Sigma_{s-1} \setminus B= (U_{s-1} \setminus B) \sqcup (\Sigma_{s-2} \setminus B)$$
involving $ \bH^{n+s}_c(\Sigma_{s-1} \setminus B; \cR)$, 
namely:
$$\cdots \to  \bH^{n+s}_c(U_{s-1} \setminus B; \cR) \to \bH^{n+s}_c(\Sigma_{s-1} \setminus B; \cR) \to \bH^{n+s}_c(\Sigma_{s-2} \setminus B; \cR) \to \cdots$$
We first show that \be\label{cl2}\bH^{n+s}_c(U_{s-1} \setminus B; \cR) \cong 0.\ee 
Indeed, the $(p,q)$-entry in the $E_2$-term of the hypercohomology spectral sequence computing $\bH^{n+s}_c(U_{s-1} \setminus B; \cR)$ is given by $$E^{p,q}_2=H^p_c(U_{s-1} \setminus B; \cH^q( \cR)),$$ and we are interested in those pairs of integers $(p,q)$ with $p+q=n+s$. 
Since a point in a $(s-1)$-dimensional stratum of $V$ has a Milnor fiber which has the homotopy type of a finite CW complex of real dimension $n-s+1$, it follows that $$\cH^q( \cR) \vert_{U_{s-1} \setminus B}\simeq 0 \ \  \text{ for any } \ q>n-s+1.$$ 
Also, by reasons of dimension, we have that $E_2^{p,q}=0$ if $p> 2s-2$. In particular, the only possibly non-trivial entries on the $E_2$-page of the above spectral sequence are those corresponding to pairs $(p,q)$ with $p\leq 2s-2$ and $q\leq n-s+1$, none of which add up to $n+s$. 
This proves \eqref{cl2}. If $s=1$, this completes the proof of \eqref{cl1} since $\Sigma_{-1}=\emptyset$. If $s>1$,  the long exact sequences for the compactly supported hypercohomology of $\cR$ associated to the disjoint union 
$$\Sigma_{\ell} \setminus B= (U_{\ell} \setminus B) \sqcup (\Sigma_{\ell-1} \setminus B),$$
$0 \leq \ell \leq s-1$, can be employed to reduce 
the proof of \eqref{cl1} to showing that
\be\label{cl3}\bH^{n+s}_c(U_{\ell} \setminus B; \cR) \cong 0\ee
for all $0 \leq \ell \leq s-1$. To prove \eqref{cl3}, we make use of the hypercohomology spectral sequence whose $E_2$-term is computed by 
$$E^{p,q}_2=H^p_c(U_{\ell} \setminus B; \cH^q( \cR)),$$ and we are interested again in those pairs of integers $(p,q)$ with $p+q=n+s$. 
Since a point in an $\ell$-dimensional stratum of $V$ has a Milnor fiber which has the homotopy type of a finite CW complex of real dimension $n-\ell$, it follows that $$\cH^q( \cR) \vert_{U_{\ell} \setminus B}\simeq 0 \ \  \text{ for any } \ q>n-\ell.$$ Moreover, by reasons of dimension, $E_2^{p,q}=0$ if $p> 2\ell$. So the only possibly non-trivial entries on the $E_2$-page are those corresponding to pairs $(p,q)$ with $p\leq 2\ell$ and $q\leq n-\ell$, none of which add up to $n+s$. This proves \eqref{cl3}, and completes the proof of \eqref{cl1} in the general case.
 
 In order to prove \eqref{rvg}, we make use of the epimorphism \eqref{nee2} as follows.
Recall that, in our notations, $U_s \setminus B$ is a disjoint union of connected strata $S_i \setminus B$ of complex dimension $s$. Each $S_i \setminus B$ has a generic transversal Milnor fiber $F_i^\pitchfork$, which has the homotopy type of a bouquet of $\mu_i^\pitchfork$ $(n-s)$-dimensional spheres. So the integral cohomology of $F_i^\pitchfork$ in concentrated in degree $n-s$. Moreover, for each $i$, there is a local system $\cL_i^\pitchfork$ on $S_i \setminus B$ with stalk $\widetilde{H}^{n-s}(F_i^\pitchfork;\bZ)$, whose monodromy is usually refered to as the {\it vertical monodromy}. This is exactly the restriction of the constructible sheaf $\cH^{n-s}( \cR)$ to 
$S_i \setminus B$. 
It then follows from the hypercohomology spectral sequence computing $\bH^{n+s}_c(U_s \setminus B; \cR) $ and by Poincar\'e duality that
\be\label{last}\bH^{n+s}_c(U_s \setminus B; \cR)  \cong \bigoplus_i \ H^{2s}_c(S_i \setminus B;\cL_i^\pitchfork) \cong \bigoplus_i \ H_0(S_i \setminus B;\cL_i^\pitchfork)\ee
which readily gives \eqref{rvg}.
$\hfill$ $\square$

\br\label{rem31}
Note that the upper bound on $b_{n+s+1}(V)$ can be formulated entirely in terms of coinvariants of vertical monodromies along the top dimensional singular strata of $V$. Indeed, if in the notations of the above proof we further let $h_i^v$ denote the vertical monodromy along $S_i \setminus B$, then each term on the right-hand side of \eqref{last} is computed by the coinvariants of $h_i^v$, i.e.,  $H_0(S_i \setminus B;\cL_i^\pitchfork) \cong \widetilde{H}^{n-s}(F_i^\pitchfork;\bZ)_{h^v_i}.$  Note that the latter statement, when combined with  \eqref{nee2}, yields an epimorphism
\be\label{nee2b} \bigoplus_i \widetilde{H}^{n-s}(F_i^\pitchfork;\bZ)_{h^v_i} \twoheadrightarrow H^{n+s}_{\varphi}(V),\ee
the summation on the left hand side being over the top dimensional singular strata of $V$. One can, moreover, proceed like in \cite{MPT} and give a more precise dependence of all (possibly non-trivial) vanishing cohomology groups $H^{k}_{\varphi}(V)$, $n \leq k \leq n+s$, in terms of the singular strata of $V$. We leave the details to the interested reader.
\er


\subsection{Isolated singularities}
Assume that the projective hypersurface $V \subset \bC P^{n+1}$ has only isolated singularities (i.e., $s=0$). 
Then the incidence variety $V_D$ is smooth since the pencil has an empty base locus, and the projection $\pi:V_D \to D$ has isolated singularities exactly at the singular points of $V$. 
 The only non-trivial vanishing homology group, $H_{n+1}^{\curlyvee}(V)$, is free, and is computed as:
\be\label{nee} H_{n+1}^{\curlyvee}(V) \cong \bigoplus_{x \in V_{\rm sing}} \widetilde{H}_{n}(F_x;\bZ) \cong \bigoplus_{x \in V_{\rm sing}} \bZ^{\mu_x},
\ee 
where $F_x$ denotes the Milnor fiber of the isolated hypersurface singularity germ $(V,x)$, with corresponding Milnor number $\mu_x$. The second isomorphism follows from the fact that $F_x$ has the homotopy type of a bouquet of $\mu_x$  $n$-spheres.

The $5$-term exact sequence \eqref{sp1} then reads as:  
\be\label{speciso}
0 \to  H_{n+1}(V_t;\bZ)  \to H_{n+1}(V;\bZ)   {\to} \bigoplus_{x \in V_{\rm sing}} \widetilde{H}_{n}(F_x;\bZ)  \overset{\alpha_{n}}{\to} H_{n}(V_t;\bZ)   \to   H_{n}(V;\bZ)   \to 0.
\ee
 Therefore Corollary \ref{corgen}(i)--(iii), together with the following bound via Theorem \ref{th2}:
 $$b_{n+1}(V) \leq 1+\sum_{x \in V_{\rm sing}} \mu_x.$$
 recover  \cite[Proposition 2.2]{ST}, which in turn is a homology counterpart of Dimca's result   \cite[Theorem 5.4.3]{Di}. 
 In fact, Dimca's result was formulated in cohomology, and it is a direct 
 consequence of the specialization sequence \eqref{spec} via Theorem \ref{th1}, together with the observation that the only non-trivial vanishing cohomology group, $H^n_{\varphi}(V)$, is  computed as:
\be\label{nee} H^n_{\varphi}(V) \cong \bigoplus_{x \in V_{\rm sing}} \widetilde{H}^n(F_x;\bZ).\ee


\br\label{rem3.3}
Let us recall here that if $V\subset \bC P^{n+1}$ is a degree $d$ reduced projective hypersurface with only isolated singularities, then its Euler characteristic is computed by the formula (e.g., see \cite[Exercise 5.3.7(i) and Corollary 5.4.4]{Di} or \cite[Proposition 10.4.2]{M}):
\be\label{chii}
\chi(V)=(n+2)-\frac{1}{d} \big[1+(-1)^{n+1}(d-1)^{n+2}\big] +(-1)^{n+1}\sum_{x \in V_{\rm sing}} \mu_x,
\ee
with $\mu_x$ denoting as before the Milnor number of the isolated hypersurface singularity germ $(V,x)$. In particular, if $V$ is a projective {\it curve} (i.e.,  $n=1$), then $H_{0}(V;\bZ) \cong \bZ$, $H_{2}(V;\bZ)\cong \bZ^r$, with $r$ denoting the number of irreducible components of $V$, and $H_{1}(V;\bZ)$ is a free group whose rank is computed from \eqref{chii} by the formula: 
\be\label{b1}b_1(V)=r+1+d^2-3d-\sum_{x \in V_{\rm sing}} \mu_x.\ee
\er


\subsection{$1$-dimensional singular locus}
  This particular case  was treated in homology in \cite[Proposition 7.7]{ST}. Let us recall the preliminaries,  in order to point out once more that in this paper we have transposed them to a fully general setting.

One starts with $V \subset \bC P^{n+1}$, a degree $d$ projective hypersurface with a singular locus $\Sigma:=V_{\rm sing}$ of complex dimension $1$. 
The singular locus $\Sigma$ consists of a union of irreducible projective curves $\Sigma_i$ and a finite set $I$ of isolated singular points. Each curve $\Sigma_i$ has a generic transversal type of transversal Milnor fiber $F_i^\pitchfork \simeq \bigvee_{\mu_i^\pitchfork} S^{n-1}$ with corresponding transversal Milnor number $\mu_i^\pitchfork$. Each $\Sigma_i$ also contains a finite set $S_i$ of special points of non-generic transversal type. One endows $V$ with the Whitney stratification whose strata are:
\begin{itemize}
\item the isolated singular points in $I$,
\item the special points in $S=\bigcup_i S_i$,
\item the (top) one-dimensional components of $\Sigma \setminus S$,
\item the open stratum $V \setminus \Sigma$.
\end{itemize}
The genericity of the pencil $\{V_{t}\}_{t\in D}$ implies that the base locus $B$ intersects each $\Sigma_i$ in a finite set $B_i$ of general points, which are not contained in $I \cup S_i$. 
The total space $V_D$ of the pencil has in this case only isolated singularities (corresponding to the points where $B$ intersects $\Sigma$), and the projection $\pi:V_D \to D$ has a $1$-dimensional singular locus $\Sigma \times \{0\}$.
  
 With the above specified landscape,  the Siersma-Tib\u ar result \cite[Proposition 7.7]{ST} reads now as the specialisation for $s=1$ of  Corollary \ref{corgenhom}, together with  the bound provided by Theorem \ref{th2}. 



\section{Examples}
In this section we work out a few specific examples. In particular, in \S\ref{quad} we show that the upper bound given by Theorem \ref{th2} is sharp, \S\ref{rathom} deals with a hypersurface which is a rational homology manifold, while \S\ref{projc} discusses the case of a projective cone on a singular curve. However, as pointed out in \cite{Di0} already in the case of isolated singularities, it is difficult  in general to compute the integral cohomology of a hypersurface by means of Corollary \ref{corgen}. It is therefore important to also develop alternative methods for exact calculations of cohomology and/or Betti numbers, e.g., see \cite{Di} for special situations.

\subsection{Singular quadrics}\label{quad} Let $n$ and $q$ be integers satisfying $4 \leq q \leq n+1$, and let 
$$f_q(x_0,\ldots x_{n+1})=\sum_{0\leq i,j \leq n+1} q_{ij} x_i x_j$$ be a quadric of rank $q:=\rk (Q)$ with $Q=(q_{ij})$. The singular locus $\Sigma$ of the quadric hypersurface $V_q=\{f_q=0\} \subset \bC P^{n+1}$ is a linear space of complex dimension $s=n+1-q$ satisfying $0 \leq s \leq n-3$. The generic transversal type for $\Sigma=\bC P^s$ is an $A_1$-singularity, so $\mu^\pitchfork=1$.
Theorem \ref{th2} yields that \be\label{ub2} b_{n+s+1}(V_q) \leq 2.\ee In what follows, we show that if the rank $q$ is even (i.e., $n+s+1$ is even), the upper bound on $ b_{n+s+1}(V_q)$ given in \eqref{ub2} is sharp. Indeed, in our notation, the quadric $V_q$ is a projective cone with vertex $\Sigma$ over a smooth quadric $W_q \subset \bC P^{n-s}$. Moreover, since $n-s\geq 3$, the homotopy version of the Lefschetz hyperplane theorem yields that $W_q$ is simply-connected (see, e.g., \cite[Theorem 1.6.5]{Di}). Let  $U=V_q \setminus \Sigma$ and consider the long exact sequence
$$ \cdots \to H^k_c(U;\bZ) \to H^k(V_q;\bZ) \to H^k(\Sigma;\bZ) \to H^{k+1}_c(U;\bZ) \to \cdots $$
Note that projecting from $\Sigma$ gives $U$ the structure of a vector bundle of rank $s+1$ over $W_q$. Let $p:U \to W_q$ denote the bundle map. Then $$H^k_c(U;\bZ)\cong H^k(W_q;Rp_!\underline{\bZ}_U)$$ can be computed by the corresponding hypercohomology spectral sequence (i.e., the compactly supported Leray-Serre spectral sequence of the map $p$), with $E^{a,b}_2=H^a(W_q;R^bp_!\underline{\bZ}_U)$. Since $\pi_1(W_q)=0$, the local system $R^bp_!\underline{\bZ}_U$ is constant on $W_q$ with stalk 
$H^b_c(\bC^{s+1};\bZ)$. Since the latter is $\bZ$ if $b=2s+2$ and $0$ otherwise, the above spectral sequence yields isomorphisms $H^k_c(U;\bZ) \cong H^{k-2-2s}(W_q;\bZ)$ if $k \geq 2s+2$ and $H^k_c(U;\bZ) \cong 0$ if $k < 2s+2$. On the other hand, $H^k(\Sigma;\bZ)=0$ if $k>2s$, so the above long exact sequence yields:
\be
H^k(V_q;\bZ)\cong 
\begin{cases}
H^k(\Sigma;\bZ) & 0 \leq k \leq 2s \\
0 &  k=2s+1 \\
H^{k-2-2s}(W_q;\bZ) & 2s+2 \leq k \leq 2n.
\end{cases}
\ee
Since $W_q$ is a smooth quadric, its integral cohomology is known from \eqref{one}, \eqref{two} and \eqref{bsm}. Altogether, this gives:
\be
H^k(V_q;\bZ)\cong 
\begin{cases}
0 & k \text{ odd} \\
\bZ & k \text{ even}, \ k \neq n+s+1 \\
\bZ^2 & k= n+s+1 \text{ even}.
\end{cases}
\ee


\subsection{One-dimensional singular locus with a two-step filtration}\label{rathom}
Let $V=\{f=0\}\subset \bC P^4$ be the $3$-fold in homogeneous coordinates $[x:y:z:t:v]$, defined by $$f=y^2z+x^3+tx^2+v^3.$$ The singular locus of $V$ is the projective line $\Sigma=\{[0:0:z:t:0] \mid z,t \in \bC\}$. By \eqref{one}, we get: $b_0(V)=1$, $b_1(V)=0$, $b_2(V)=1$. Since $V$ is irreducible, \eqref{top} yields: $b_6(V)=1$. 
We are therefore interested to understand the Betti numbers $b_3(V)$, $b_4(V)$ and $b_5(V)$.

It was shown in \cite[Example 6.1]{M0} that $V$ has a Whitney stratification with strata: $$S_3:=V \setminus \Sigma, \ \ S_1:=\Sigma \setminus [0:0:0:1:0], \ \ S_0:=[0:0:0:1:0],$$ giving $V$ a two-step filtration $V \supset \Sigma \supset [0:0:0:1:0].$

The transversal singularity for the top singular stratum $S_1$ is the Brieskorn type singularity $y^2+x^3+v^3=0$ at the origin of $\bC^3$ (in a normal slice to $S_1$), with corresponding transversal Milnor number $\mu_1^\pitchfork =4$. So Theorem \ref{th2} yields that $b_5(V) \leq 5$, while Corollary \ref{corgen} gives $b_3(V) \leq 10$. As we will indicate below, the actual values of $b_3(V)$ and $b_5(V)$ are zero.

It was shown in \cite[Example 6.1]{M0} that the hypersurface $V$ is in fact a $\bQ$-homology manifold, so it satisfies Poincar\'e duality over the rationals. In particular, $b_5(V)=b_1(V)=0$ and $b_4(V)=b_2(V)=1$. To determine $b_3(V)$, it suffices to compute  the Euler characteristic of $V$, since $\chi(V)=4-b_3(V)$. Let us denote by $Y\subset \bC P^4$ a smooth $3$-fold which intersects the Whitney stratification of $V$ transversally. Then \eqref{chi} yields that $\chi(Y)=-6$ and we have by \cite[(10.40)]{M} that
\be\label{plug}
\chi(V)=\chi(Y)-\chi(S_1 \setminus Y) \cdot \mu_1^\pitchfork -\chi(S_0) \cdot (\chi(F_0)-1),
\ee
where $F_0$ denotes the Milnor fiber of $V$ at the singular point $S_0$. As shown in \cite[Example 6.1]{M0}, $F_0 \simeq S^3 \vee S^3$. So, using the fact that the general $3$-fold $Y$ intersects $S_1$ at $3$ points, we get from \eqref{plug} that $\chi(V)=4$. Therefore, $b_3(V)=0$, as claimed. Moreover, since $H^3(V;\bZ)$ is free, this also shows that in fact $H^3(V;\bZ)\cong 0$.

\br Note that the hypersurface of the previous example has the same Betti numbers as $\bC P^3$. This fact can also be checked directly, by noting that the monodromy operator acting on the reduced homology of the Milnor fiber of $f$ at the origin in $\bC^5$ has no eigenvalue equal to $1$ (see \cite[Corollary 5.2.22]{Di}).

More generally, consider a degree $d$ homogeneous polynomial $g(x_0,\ldots, x_n)$ with associated Milnor $F_g$ such that the monodromy operator $h_*$ acting on $\widetilde{H}_*(F_g;\bQ)$ is the identity. Then the hypersurface $V=\{g(x_0,\ldots, x_n)+x_{n+1}^{d}=0\}\subset \bC P^{n+1}$ has the same $\bQ$-(co)homology as $\bC P^n$. For example, the hypersurface $V_n=\{x_0x_1\ldots x_n+x_{n+1}^{n+1}=0\}$ has singularities in codimension $2$, but the same $\bQ$-(co)homology as $\bC P^n$. However, $V_n$ does not have in general the $\bZ$-(co)homology of $\bC P^n$; indeed, $H^3(V_2;\bZ)$ contains  $3$-torsion (cf. \cite[Proposition 5.4.8]{Di}).
\er

\subsection{Projective cone on a curve}\label{projc}
The  projective curve $C=\{xyz=0\}\subset \bC P^2$ has three irreducible components and three singularities of type $A_1$ (each having a corresponding Milnor number equal to $1$). Therefore, by Remark \ref{rem3.3} and formula \eqref{b1}, the integral cohomology of $C$ is given by:
$$H^0(C;\bZ)\cong \bZ, \ H^1(C;\bZ) \cong \bZ, \ H^2(C;\bZ)\cong \bZ^3.$$
The projective cone on $C$ is the surface $V=\{xyz=0\}\subset \bC P^3$. The singular locus of $V$ consists of three projective lines intersecting at the point $[0:0:0:1]$, each having a (generic) transversal singularity of type $A_1$, i.e., with corresponding transversal Milnor number equal to $1$. 
By \cite[(5.4.18)]{Di}, we have that  $$H^k(V;\bZ) \cong H^{k-2}(C;\bZ), \ \ \text{for all} \ k \geq 2.$$
Together with \eqref{one}, this yields:
\be\label{comp1}
H^0(V;\bZ)\cong \bZ, \ H^1(V;\bZ)\cong 0, \ H^2(V;\bZ)\cong \bZ, \ H^3(V;\bZ) \cong \bZ, \ H^4(V;\bZ)\cong \bZ^3.
\ee
By Theorem \ref{th1}, the only non-trivial  vanishing cohomology groups of $V$ are $H^2_{\varphi}(V)$, which is free, and $H^3_{\varphi}(V)$. These can be explicitly computed by using 
\eqref{bsm}, \eqref{sp1} and \eqref{comp1}, to get: $$H^2_{\varphi}(V)\cong \bZ^7, \ H^3_{\varphi}(V)\cong \bZ^2$$ (compare with \cite[Example 7.5]{ST}).


\section{Supplement to the Lefschetz hyperplane theorem and applications}\label{supLHT}
In this section, we give a new proof of Kato's result mentioned in the Introduction. Our proof is different from that of \cite[Theorem 5.2.11]{Di}, and it relies on a supplement to the Lefschetz hyperplane section theorem (Theorem \ref{thapi}), which is proved in Theorem \ref{thap} below.

\subsection{A supplement to the Lefschetz hyperplane theorem}
In this section, we prove the following result of Lefschetz type:
\bt\label{thap}
Let $V \subset \bC P^{n+1}$ be a reduced complex projective hypersurface with $s=\dim V_{\rm sing}$ the complex dimension of its singular locus. (By convention, we set $s=-1$ if $V$ is nonsingular.) Let $H \subset \bC P^{n+1}$ be a generic hyperplane (i.e., transversal to a Whitney stratification of $V$), and denote by $V_H:=V\cap H$ the corresponding hyperplane section of $V$.
Then
\be\label{34ap}
H^k(V,V_H; \bZ)=0 \ \ \text{for} \ \ k < n \ \ \text{and} \ \ n+s+1 < k < 2n.
\ee
Moreover, $H^{2n}(V,V_H; \bZ)\cong\bZ^r$, where $r$ is the number of irreducible components of $V$, and $H^{n}(V,V_H; \bZ)$ is (torsion-)free.
\et
\begin{proof}
Let us first note that the long exact sequence for the cohomology of the pair $(V,V_H)$ together with \eqref{top} yield that:
$$H^{2n}(V,V_H; \bZ)\cong H^{2n}(V;\bZ)\cong\bZ^r.$$
Moreover, we have isomorphisms:
$$H^k(V,V_H; \bZ) \cong H^k_c(V^a;\bZ),$$
where $V^a:=V\setminus V_H$. 
Therefore, the vanishing in \eqref{34ap} for $k<n$ is a consequence of the Artin vanishing theorem (e.g., see \cite[Corollary 6.0.4]{Sc}) for the perverse sheaf $\underline{\bZ}_{V^a}[n]$ (cf. \cite{Le}) on the affine hypersurface $V^a$ obtained from $V$ by removing the hyperplane section $V_H$. Indeed,
$$H^k_c(V^a;\bZ)=\bH^{k-n}_c(V^a;\underline{\bZ}_{V^a}[n]) \cong 0$$ for all $k-n<0$. (Note that vanishing in this range is equivalent to the classical Lefschetz hyperplane section theorem.)

Since $V$ is reduced, we have that $s<n$. If $n=s+1$ then $n+s+1=2n$ and there is nothing else to prove in \eqref{34ap}. So let us now assume that $n>s+1$.  For $n+s+1<k<2n$, we have the following sequence of isomorphisms:
\be\label{35} \begin{split}
H^k(V, V_H; \bZ) &\cong  H^k(V \cup H, H; \bZ) \\ 
&\cong H_{2n+2-k} (\bC P^{n+1}\setminus H, \bC P^{n+1} \setminus (V \cup H); \bZ) \\
&\cong H_{2n+1-k}(\bC P^{n+1}\setminus (V \cup H); \bZ),
\end{split}
\ee
where the first isomorphism follows by excision, the second is an application of the Poincar\'e-Alexander-Lefschetz duality, and the third follows from the cohomology long exact sequence of a pair.
Set 
$$U=\bC P^{n+1}\setminus (V \cup H),$$ and let $L = \bC P^{n-s}$ be a generic linear subspace (i.e., transversal to both $V$ and $H$). Then, by transversality, $L \cap V$ is a nonsingular hypersurface in $L$, transversal to the hyperplane at infinity $L \cap H$ in $L$. Therefore, $U \cap L=L\setminus (V \cup H) \cap L$ has the homotopy type of a wedge
$$U \cap L \simeq S^1 \vee S^{n-s} \vee \ldots \vee S^{n-s},$$
e.g., see \cite[Corollary 1.2]{Lib}. Thus, by the Lefschetz hyperplane section theorem (applied $s+1$ times), we obtain:
$$H_i(U;\bZ) \cong H_i(U \cap L;\bZ) \cong 0$$
for all integers $i$ in the range $1 < i < n-s$. Substituting $i = 2n+1-k$ in \eqref{35}, we get that
$H^k(V, V_H; \bZ)\cong 0$ for all integers $k$ in the range $n+s+1 < k < 2n$. 

It remains to show that $H^{n}(V,V_H; \bZ)\cong H^n_c(V^a;\bZ)\cong \bH^{0}_c(V^a;\underline{\bZ}_{V^a}[n])$ is (torsion-)free. This follows as in the proof of Theorem \ref{th1} since the affine hypersurface $V^a$ has rectified homological depth equal to its complex dimension $n$.
This completes the proof of the theorem.
\end{proof}

Theorem \ref{thap} and the Universal Coefficient Theorem now yield the following consequence:
\bc\label{corap}
In the notations of Theorem \ref{thap} we have that:
\be\label{36ap}
H_k(V,V_H; \bZ)=0 \ \ \text{for} \ \ k < n \ \ \text{and} \ \ n+s+1 < k < 2n.
\ee
Moreover, $H_{2n}(V,V_H; \bZ)\cong\bZ^r$, where $r$ is the number of irreducible components of $V$.
\ec


\subsection{Kato's theorem for hypersurfaces}
The isomorphism \eqref{two} from the introduction was originally proved by Kato \cite{Ka}, and it holds more generally for complete intersections. We derive it here as a consequence of Theorem \ref{thap}.

\bt[Kato]\label{Kato}
Let $V\subset \bC P^{n+1}$ be a reduced degree $d$ complex projective hypersurface with $s=\dim V_{\rm sing}$ the complex dimension of its singular locus. (By convention, we set $s=-1$ if $V$ is nonsingular.) Then
\be\label{twoap}
H^k(V;\bZ) \cong H^k( \bC P^{n+1};\bZ) \ \  \text{for all} \ \  n+s+2\leq k\leq 2n.
\ee
Moreover, if $j:V \hookrightarrow \bC P^{n+1}$ denotes the inclusion, the induced cohomology homomorphisms 
\be\label{threeap}
j^k:H^k( \bC P^{n+1};\bZ) \lra H^k(V;\bZ), \ \ n+s+2\leq k\leq 2n,
\ee
are given by multiplication by $d$ if $k$ is even.
\et

\begin{proof}
The statement of the theorem is valid only if $n\geq s+2$, so in particular we can assume that $V$ is irreducible and hence $H^{2n}(V;\bZ)\cong \bZ$.  
Moreover, the fact that $j^{2n}$ is multiplication by $d=\deg(V)$ is true  regardless of the dimension of singular locus, see \cite[(5.2.10)]{Di}.
If $n=s+2$ there is nothing else to prove, so we may assume  (without any loss of generality) that $n \geq s+3$. 

We next proceed by induction on $s$. 

If $V$ is nonsingular (i.e., $s=-1$), the assertions are well-known for any $n \geq 1$. We include here a proof for completeness. The isomorphism \eqref{twoap} can be obtained in this case from the Lefschetz isomorphism \eqref{one}, its homology analogue, and Poincar\'e duality.  The statement about $j^k$ can also be deduced from \eqref{one} and Poincar\'e duality, but we include here a different argument inspired by \cite{Di}. Consider the isolated singularity at the origin for the affine cone $CV \subset \bC^{n+2}$ on $V$, and the corresponding link $L_V:=S^{2n+3} \cap CV$, for $S^{2n+3}$ a small enough sphere at the origin in $\bC^{n+2}$. Then $L_V$ is a $(n-1)$-connected closed oriented manifold of real dimension $2n+1$, so its only possibly nontrivial integral (co)homology appears in degrees $0$, $n$, $n+1$ and $2n+1$. The Hopf fibration $S^1 \hookrightarrow S^{2n+3} \lra \bC P^{n+1}$ induces by restriction to $CV$ a corresponding Hopf fibration for $V$, namely $S^1 \hookrightarrow L_V \lra V$. Then for any $n+1 \leq k \leq 2n-2$, the cohomology Gysin sequences for the diagram of fibrations 
$$\xymatrix{
S^{2n+3} \ar[r]  & \bC P^{n+1}\\
L_V \ar[u] \ar[r] & V \ar[u].
}$$
yield commutative diagrams (with $\bZ$-coefficients):
\be\label{Gys} \CD
0=H^{k+1}(S^{2n+3}) @>>> H^k(\bC P^{n+1}) @>{\psi}>{\cong}>  H^{k+2}(\bC P^{n+1}) @>>> H^{k+2}(S^{2n+3})=0\\
@VVV @V {j^k}VV  @V {j^{k+2}}VV @VVV \\
0=H^{k+1}(L_V) @>>> H^k(V) @>{\psi_V}>{\cong}>  H^{k+2}(V) @>>> H^{k+2}(L_V)=0\\
\endCD \ee
Here, if $k=2\ell$ is even, the isomorphism $\psi$ is the cup product with the cohomology generator $a\in H^2(\bC P^{n+1};\bZ)$, and similarly, $\psi_V$ is the cup product with $j^2(a)$. The assertion about $j^k$ follows now from \eqref{Gys} by decreasing induction on $\ell$, using the fact mentioned at the beginning of the proof that $j^{2n}$ is given by multiplication by $d$.

Let us next choose a generic hyperplane $H \subset \bC P^{n+1}$ (i.e., $H$ is transversal to a Whitney stratification of $V$), and set as before $V_H=V \cap H$. It then follows from Theorem \ref{thap} and the cohomology long exact sequence of the pair $(V,V_H)$ that $H^{2n-1}(V;\bZ) \cong 0$. It therefore remains to prove \eqref{twoap} and the corresponding assertion about $j^k$ for $k$ in the range for $n+s+2\leq k\leq 2n-2$.
Let us consider the commuting square
$$\CD
V_H @> {\delta} >> H=\bC P^n\\
@V {\gamma}VV @VVV\\\
V @>> {j} > \bC P^{n+1}
\endCD$$
and the induced commutative diagram in cohomology:
\be\label{di} \CD
H^k(\bC P^{n+1};\bZ) @> {j^k} >> H^k(V;\bZ) \\
@V {\cong}VV @VV{\gamma^k}V\\\
H^k(\bC P^{n};\bZ) @>> {\delta^k} > H^k(V_H;\bZ)
\endCD \ee
By Theorem \ref{thap} and the cohomology long exact sequence of the pair $(V,V_H)$ we get that $\gamma^k$ is an isomorphism for all integers $k$ in the range $n+s+2\leq k\leq 2n-2$. Moreover, since $V_H \subset \bC P^n$ is a degree $d$ reduced projective hypersurface with a $(s-1)$-dimensional singular locus (by transversality), the induction hypothesis yields that $H^k(V_H;\bZ) \cong H^k(\bC  P^n;\bZ)$ for $n+s \leq k \leq 2n-2$ and that, in the same range and for $k$ even, the homomorphism 
$\delta^k$ is given by multiplication by $d$. The commutativity of the above diagram \eqref{di} then yields \eqref{twoap} for all integers $k$ satisfying $n+s+2\leq k\leq 2n-2$, and the corresponding assertion about the induced homomorphism $j^k$ for $k$ even in the same range. This completes the proof of the theorem.
\end{proof}

\br
Let us remark here that the proof of Kato's theorem in \cite[Theorem 5.2.11]{Di} relies on the Kato-Matsumoto result \cite{KM} on the connectivity of the Milnor fiber of the singularity at the origin of the affine cone $CV \subset \bC^{n+2}$. 
\er

\br\label{Katoh}
 One can prove the homological version of Theorem \ref{Kato} in the similar manner,  namely by using Corollary \ref{corap} instead of Theorem \ref{thap}. This yields the isomorphisms: 
\be
H_k(V;\bZ)  \cong H_k( \bC P^{n+1};\bZ)  \ \  \text{for all} \ \  n+s+2\leq k\leq 2n,
\ee
and the homomorphisms induced  by the inclusion $j:V\hookrightarrow \bC P^{n+1}$ in homology are given in this range (and for $k$ even) by multiplication by $d=\deg(V)$.
\er

\br
We already noted that Theorem \ref{th1} yields the isomorphism \eqref{two} of Kato's theorem (see Corollary \ref{corgen}(i)). On the other hand, Kato's Theorem \ref{Kato} may be used to obtain a weaker version of Theorem \ref{th1} by more elementary means. Indeed, in the notations from the Introduction consider the diagram:
$$ \CD
H^k(\bC P^{n+1};\bZ) @> {\cong}>> H^k(\bC P^{n+1} \times D;\bZ) @> {b^k} >> H^k(V_D;\bZ)  @> {c^k} >> H^k(V_t;\bZ)\\
                           @.       @.   @V {\cong}VV @. \\
@.       @.  H^k(V;\bZ)
\endCD 
$$
and let $a^k:=c^k \circ b^k$. By Theorem \ref{Kato}, we have that:
\begin{itemize}
\item[(i)] $a^k$ is the multiplication by $d$ if $k>n$ even and an isomorphism for $k<n$;
\item[(ii)] $b^k$ is the multiplication by $d$ if $n+s+2\leq k \leq 2n$ ($k$ even) and an isomorphism for $k<n$.
\end{itemize}
Therefore, $c^k$ is an isomorphism if $n+s+2 \leq k \leq 2n$ or $k<n$. The cohomology long exact sequence of the pair $(V_D,V_t)$ then yields that $H^k_{\varphi}(V)\cong H^{k+1}(V_D,V_t;\bZ)\cong 0$ for all integers $k \notin [n-1,n+s+1]$.
\er


\begin{thebibliography}{99}

\bibitem{Di0} Dimca, Alexandru, {\it On the homology and cohomology of complete intersections with isolated singularities} Compositio Math. 58 (1986), no. 3, 321--339.

\bibitem{Di} Dimca, Alexandru, {\it Singularities and Topology of Hypersurfaces}, Universitext, Springer, 1992.

\bibitem{Di1} Dimca, Alexandru, {\it Sheaves in Topology}, Universitext, Springer-Verlag, Berlin, 2004.

\bibitem{Ka} Kato, Mitsuyoshi, {\it Topology of $k$-regular spaces and algebraic sets},  Manifolds -- Tokyo 1973 (Proc. Internat. Conf. on Manifolds and Related Topics in Topology), pp. 153--159. Univ. Tokyo Press, Tokyo, 1975.

\bibitem{KM} Kato, Mitsuyoshi, Matsumoto, Yukio,
{\it On the connectivity of the Milnor fiber of a holomorphic function at a critical point}, Manifolds--Tokyo 1973 (Proc. Internat. Conf., Tokyo, 1973), pp. 131--136. Univ. Tokyo Press, Tokyo, 1975.

\bibitem{Le} L\^e, D\~ung Tr\'ang, {\it Sur les cycles \'evanouissants des espaces analytiques}, C. R. Acad. Sci. Paris S\'er. A-B 288(4), A283--A285 (1979).

\bibitem{Lib} Libgober, Anatoly, {\it Homotopy groups of the complements to singular hypersurfaces, II}, Ann. of Math. (2) 139 (1994), 117--144.


\bibitem{Ma} Massey, David, {\it  Natural commuting of vanishing cycles and the Verdier dual}, Pacific J. Math. 284 (2016), no. 2, 431--437.

\bibitem{M0} Maxim, Laurentiu, {\it Intersection homology and Alexander modules of hypersurface complements}, Comment. Math. Helv. 81 (2006), no. 1, 123--155.

\bibitem{M} Maxim, Laurentiu, {\it Intersection Homology \& Perverse Sheaves, with Applications to Singularities}, Graduate Texts in Mathematics, Vol. 281, Springer, 2019.

\bibitem{MPT} Maxim, Laurentiu, Paunescu, Laurentiu, Tibar, Mihai, {\it The vanishing cohomology of non-isolated hypersurface singularities},
arXiv:2007.07064

\bibitem{MSS} Maxim, Laurentiu, Saito, Morihiko, Sch\"urmann, J\"org, {\it Hirzebruch-Milnor classes of complete intersections}, Adv. Math. 241 (2013), 220--245.

\bibitem{Mi} Miller, John L., {\it Homology of complex projective hypersurfaces with isolated singularities}, Proc. Amer. Math. Soc. 56 (1976), 310--312.

\bibitem{PP} Parusi\'nski, Adam, Pragacz, Piotr, 
{\it Characteristic classes of hypersurfaces and characteristic cycles}, 
J. Algebraic Geom. 10 (2001), no. 1, 63--79.

\bibitem{ST0} Siersma, Dirk, Tib\u{a}r, Mihai,  {\it Milnor fibre homology via deformation}, Singularities and computer algebra, 305--322, Springer, Cham, 2017. 

\bibitem{ST} Siersma, Dirk, Tib\u{a}r, Mihai,  {\it Vanishing homology of projective hypersurfaces with $1$-dimensional singularities}, Eur. J. Math. 3 (2017), no. 3, 565--586.

\bibitem{Sc} Sch\"urmann, J\"org, {\it Topology of Singular Spaces and Constructible Sheaves}, Birkh\"auser, Monografie Matematyczne 63, 2003.

\end{thebibliography}
\end{document}